\documentclass[11pt]{article}
\usepackage{amssymb}
\usepackage{latexsym}

\newtheorem{thm}{Th\'eor\`eme}[section]

\newtheorem{lemma}[thm]{Lemme}

\begin{document}

\title{\bf Le rang des tissus de Nakai. } 

\author{ Jean-Paul Dufour et Daniel Lehmann \\}

\maketitle

\begin{abstract}
According to Alain H\'enaut, a planar 4-web is called Nakai's web if the cross-ratio of the tangents to the four foliations at each point is constant and if it has no hexagonal 3-subweb. We prove that Nakai's webs have rank 0 or 1. We give examples with rank 1 and present a universal way to build such examples.

\end{abstract}

\section {\bf Introduction. }

Nous travaillons ici  dans le cadre analytique r\'eel  aussi bien que  dans le cadre holomorphe. Nous travaillerons toujours "localement", c'est \`a dire au voisinage d'un point du plan que nous appellerons "l'origine" et qui sera le point de coordonn\'ees $(0,0)$ dans tous les syst\`emes de coordonn\'ees que nous utiliserons.

 Nous \'etudions les 4-tissus du plan dont le "bi-rapport" (i.e. le bi-rapport des tangentes aux feuilles) en chaque point est constant donc ind\'ependant du point.

En conformit\'e avec A. H\'enaut, nous appelons {\bf tissu de Nakai} un 4-tissu plan dont le bi-rapport est constant et dont aucun des 3-sous-tissus n'est hexagonal.

Rappelons que Isao Nakai a \'etudi\'e ce type de tissu dans [N]. Il est facile de voir que, pour un 4-tissu du plan dont le bi-rapport  est constant, il suffit que l'un des 3-sous-tissus soit non hexagonal pour que tous  soient non-hexagonaux et, de plus, ils ont tous la m\^eme courbure de Blaschke. Un des r\'esultats non triviaux de [N] est que, si un 4-tissu est tel que tous ses 3-sous-tissus ont la m\^eme courbure de Blaschke non-nulle, alors le tissu a un bi-rapport constant.

Dans ce texte, un 4-tissu planaire sera  donn\'e par  quatre fonctions $f_1,$ $f_2,$ $f_3,$ et $f_4,$ qui sont int\'egrales premi\`eres de ses feuilletages (et sera not\'e $(f_1,f_2,f_3,f_4)$) ou par quatre 1-formes $\omega_1,$ $\omega_2,$ $\omega_3,$ et $\omega_4$ dont les noyaux sont  tangents aux feuilles (et sera not\'e $(\omega_1,\omega_2,\omega_3,\omega_4)).$ 
  
Rappelons que les relations ab\'eliennes d'un 4-tissu planaire  $(f_1,f_2,f_3,f_4)$ sont les familles de fonctions $a_1,a_2,a_3,a_4$ 
d'une variable (d\'efinies \`a constante additive pr\`es) telles que
$$a_1(f_1)+a_2(f_2)+a_3(f_3)+a_4(f_4)\equiv cte.$$
Il est \'equivalent de dire que ce sont les familles $ (\omega_1,\omega_ 2, \omega_3,\omega_4)$
de 1-formes ferm\'ees,
respectivement colin\'eaires \`a $df_1,df_2, df_3, df_4,$ et dont la somme est identiquement nulle : elles forment donc
un espace vectoriel, dont la dimension s'appelle le {\sl rang du tissu}.

Le rang maximal d'un 4-tissu planaire est 3, mais nous savons d\'ej\`a qu'un tissu de Nakai ne peut
pas \^etre de rang 3, puisque sa courbure de Blaschke g\'en\'eralis\'ee,
c'est-\`a-dire la somme des courbures de Blaschke des quatre 3-sous-tissus, qui est aussi \'egale \`a la trace de la courbure d\'efinie par Pantazi [Pa], H\'enaut [H] ou Pirio [Pi], ne peut pas s'annuler. Dans cet article, nous montrerons que ce rang est au plus 1, qu'il peut \^etre \'egal \`a 1, mais qu'il est g\'en\'eriquement nul et nous construirons des exemples de tissus de Nakai de rang 1.

La technique d\'evelopp\'ee dans [DL] permet de voir que le rang des tissus de Nakai est g\'en\'eriquement 0 (on y construit une matrice $20\times 20$ dont le d\'eterminant est en g\'en\'eral non-nul, interdisant l'existence de relations ab\'eliennes non triviales). En utilisant ces techniques on retrouve que le rang est au plus 2 mais les calculs deviennent trop lourds pour aller plus loin. Nous d\'emontrerons nos r\'esultats par une m\'ethode sp\'ecifique.

Nous remercions A. H\'enaut de nous avoir pos\'e le probl\`eme.

\section{\bf Relations ab\'eliennes des 4-tissus \`a bi-rapport constant.}

A isomorphisme pr\`es, tout 4-tissu plan de bi-rapport constant  s'obtient comme suit : 

On part d'un 3-tissu $W= (dx,dy,dy-p.dx),$ o\`u $p=p(x,y)$ est une fonction non-nulle \`a l'origine  arbitraire, et on lui rajoute le feuilletage engendr\'e par la quatri\`eme forme $dy-C.p.dx$ o\`u $C$ est un nombre diff\'erent de 0 et 1. On notera $W_C$ ce 4-tissu.

Dans cette section nous ne faisons aucune hypoth\`ese concernant les courbures des 3-sous-tissus.  donc $W_C$ n'est pas forc\'ement un tissu de Nakai.

Les relations ab\'eliennes de $W_C$ sont obtenues \`a partir de quatre fonctions $r,$ $s,$ $c$ et $e$  telles que 

$$r.dx+s.dy+c.(dy-p.dx)+e.(dy-C.p.dx)=0,$$ 
$$d(r.dx)=0,\ \ \ \ d(s.dy)=0,\ \ \ \ d(c.(dy-p.dx))=0,\ \ \ \ d(e.(dy-C.p.dx))=0.$$

Ce syst\`eme d'\'equations se r\'e\'ecrit comme suit.
Les fonctions $r$ et $s$ sont des fonctions d'une seule variable, respectivement $x$ et $y$ et on a 
 $$c'_x-(c.p)'_y=0,\ \ \ \  e'_x-C.(e.p)'_y =0,\ \ \ \  s+c+e=0,\ \ \ \  r-p.c-C.p.e=0,$$
(pour toute fonction $g$, la notation $g'_x$ (resp. $g'_y$) d\'esigne la d\'eriv\'ee de $g$ par rapport \`a $x$ (resp. $y$)).

Les deux derni\`eres relations permettent d'exprimer $c$ et $e$ en en fonction de $r$, $s$ et $p.$  On trouve en particulier
$$c=(C.s+r/p)/(1-C).$$ Alors l'\'equation $c'_x-(c.p)'_y=0$ devient 

$$(*)\ \ \ \ r'(x).p-r(x).p'_x+C.(s'(y).p^3+s(y).p^2.p'_y)=0.$$

On v\'erifie que la relation $ e'_x-C.(e.p)'_y =0$ est aussi \'equivalente \`a $(*)$.

En conclusion de ces calculs, nous retiendrons que les relations ab\'eliennes de notre tissu $W_C$ sont donn\'ees par les solutions $(r,s)$ de $(*).$ 
\begin{thm} Soient deux nombres $C$ et $D$ non-nuls et non \'egaux \`a 1. Il y a une bijection entre les relations ab\'eliennes de $W_C$ et celles de $W_D$\end{thm}

La preuve, \'evidente, est bas\'ee sur le fait que si $(r,s)$ est solution de $(*)$ alors $(r,D.s/C)$ est solution de l'\'equation obtenue en rempla\c cant $C$ par $D$ dans $(*).$ 

Ceci justifie le fait que, dans la suite de ce texte, on ne s'int\'eresse qu'aux 4-tissus de bi-rapport constant \'egal \`a -1, que nous appellerons "harmoniques." Quand nous parlerons de la relation $(*)$ il sera sous-entendu que nous imposerons $C=-1.$

\section{Rang maximum des tissus de Nakai.}

Dans cette section nous nous int\'eressons aux tissus de Nakai. Ce sont les 4-tissus du type de ceux de la section pr\'ec\'edente, mais l'on impose que $W$ ne soit pas hexagonal. De fa\c con plus pr\'ecise nous imposerons que la courbure de Blaschke de $W$  ne soit pas nulle \`a l'origine.

Rappelons le r\'esultat classique suivant.

\begin{lemma} On consid\`ere un 3-tissu  ayant  une courbure de Blaschke non-nulle \`a l'origine. Il est isomorphe \`a un unique  tissu $W=(dx,dy,dy-p.dx)$ avec $p$ de la forme $$p(x,y)=1+xy(1+h(x,y))$$ o\`u $h$ est une fonction qui s'annule \`a l'origine.\end{lemma}

On d\'emontre ce lemme en faisant des changements de variables du type $(x,y)\mapsto (X(x),Y(y)).$ La forme qu'il donne \`a $p$ est appel\'ee la "forme normale".

Pour trouver le rang des tissus de Nakai il faut donc \'etudier l'\'equation $(*)$ dans le cas o\`u $p$ a cette forme normale. Nous allons voir qu'il suffit d'\'etudier le d\'eveloppement de Taylor d'ordre 3 du premier membre de $(*)$. 

Pour cela on note $ r_i$ (resp. $s_i$) les coefficients de Taylor de $r$ (resp. $s$). On suppose aussi que le d\'eveloppement de Taylor d'ordre 2 de la fonction $h$ qui apparait dans la forme normale de $p$ est
 $$ax+by+ux^2+vxy+wy^2.$$

 Un calcul \'el\'ementaire  donne les r\'esultats suivants.

Les termes constants donnent $s_1=r_1.$ Les termes de degr\'e 1 donnent
$$s_2 = -(1/2)r_0,\ \ \ \ r_2 = (1/2)s_0.$$
Les termes de degr\'e 2 donnent
$$r_1 = -(1/2)b.s_0-(1/2)a.r_0, \ \ \   s_3 = -(1/3)b.r_0,\ \ \ \  r_3 = (1/3)a.s_0.$$
Les termes de degr\'e 3 donnent, d'une part,
$$r_4 = (1/4)u.s_0,\ \ \ \  s_4 = -(1/4)w.r_0,$$
d'autre part les deux relations 
$$(-5a^2+6u)r_0+(3-5b.a+4v)s_0=0,\ \ \ (-5b.a+4v-7)r_0+(6w-5b^2)s_0=0.$$

Cela donne un syst\`eme lin\'eaire de variables $r_0$ et $s_0.$ Comme la matrice de ce syst\`eme est de rang au moins 1, et en g\'en\'eral de rang 2. On voit qu'en g\'en\'eral on a $r_0=s_0=0$ et un raisonnement par r\'ecurrence permet de voir que tous les $r_i$ et $s_i$ sont nuls, donc que $s$ et $r$ sont nuls, donc que la relation ab\'elienne cherch\'ee est triviale. Dans le cas o\`u notre matrice est de rang 1 cela montre que l'on a soit $s_0=\lambda. r_0$ soit $r_0=\mu. s_0$ et le m\^eme raisonnement par r\'ecurrence montrent que les $r_i$ et les $s_i$ ont tous le facteur $r_0$ ou le facteur $s_0$. Cela montre que les relations ab\'eliennes sont proportionelles. D'o\`u le th\'eor\`eme suivant.

\begin{thm} Les tissus de Nakai sont en g\'en\'eral de rang 0 et, au maximum, de rang 1.\end{thm}

{\bf Remarque}. Les calculs pr\'ec\'edents prouvent aussi que, si un tissu de Nakai a une relation ab\'elienne non triviale elle provient forc\'ement d'une \'equation $(*)$ avec $r_0$ ou $s_0$ non nul.

 \section{Tissus de  Nakai de rang un.}
    \subsection{Un premier exemple (harmonique)}

   Pour toute      fonction  $u$ d'une variable, d\'efinissons le 4-tissu $W_u=(x,\ y,\ f,\ g)$, avec    
   $$f(x,y):= \frac{x^3}{6}+y+x.u(xy)\hbox{ \ et\ \ }g(x,y):= \frac{x^3}{6}+y-x.u(xy).$$
  Ce tissu poss\`ede une relation ab\'elienne \'evidente non triviale $$\frac{x^3}{3}+ 2y-f-g\equiv 0 $$. 

 Nous allons calculer    $u$ 
   de fa\c con que  $W_u$ soit un tissu harmonique. Pour cel\`a \'ecrivons d'abord que le bi-rapport  $(0,\infty,-\frac{f'_x}{f'_y},-\frac{g'_x}{g'_y})$ est \'egal \`a -1 :
   $$f'_x.g'_y+g'_x.f'_y=0. $$
 On obtient  :$$x^2\Bigl(
   u'(xy).\bigl(u(xy)+xy.u'(xy)\bigr)-\frac{1}{2}\Bigr)=0.$$
  Posant  $t=xy$, la  fonction $u$ doit \^etre solution de    l'\'equation diff\'erentielle  $t. (\frac{du}{dt})^2+u \frac{du}{dt}-\frac{1}{2}=0. $ On remarque que les \'eventuelles solutions de cette \'equation v\'erifient $u(0).u'(0)=\frac{1}{2}$ ;  en particulier $u'(0)$ n'est pas nul (il existe donc un nombre $a$ non nul tel que $u(0)=a$ et $u'(0)=\frac{1}{2a}$), et la fonction $t\mapsto u(t)$ doit admettre une fonction r\'ecipoque $u\mapsto t(u)$ solution de l'\'equation diff\'erentielle
 $$\Bigl(\frac{dt}{du}\Bigr)^2-2u. \frac{dt}{du} -2t=0.$$ En dehors des  points $(t,u) $ tels que  $u^2+2t=0$, l'\'equation pr\'ec\'edente 
se scinde en deux  \'equations diff\'erentielles$$\frac{dt}{du}=u\pm (u^2+2t) ^{1/2}, $$ dont on cherche les solutions avec condition initiale $t(a)=0$. L'une d'elles admet la solution $t=0$ qui ne convient pas ;  l'autre admet, pour tout $a\neq 0$ choisi arbitrairement,   une fonction    solution $u\mapsto t(u)$ v\'erifiant  $t(a)=0$ ; et cette solution   admet bien une fonction r\'eciproque au voisinage de $a$. 

Il reste \ a montrer que $W_u$ n'a pas de   3-sous-tissu hexagonal. On v\'erifie ais\'ement que le     3-sous-tissu   $(x,y,f)$ admet \`a l'origine la courbure $\frac{1}{a^2}\ dx\wedge dy$, et n'est donc pas hexagonal ;  les  autres ne le sont pas non plus,  en vertu de la remarque faite dans l'introduction.

On peut aussi remarquer qu'un isomorphisme $(x,y)\mapsto (\lambda x,y/\lambda)$ permet toujours de nous ramener au cas o\`u $a$ est 1.
 
\subsection{M\'ethode g\'en\'erale pour construire les tissus de Nakai de rang 1.}

    Notre outil sera encore  l'\'equation (*). Mais, alors que dans la d\'emonstration de th\'eor\`eme 3.2 nous fixions $p$ et cherchions $r$ et $s$, ici nous faisons le contraire.  Nous fixons $r$ et $s$ et cherchons un $p$ qui pourrait nous donner un exemple de tissu de Nakai. La remarque faite \`a la suite de l'\'enonc\'e du th\'eor\`eme 3.2 sera ici fondamentale : elle nous dit que l'on doit se restreindre \`a des fonctions $r$ et $s$ qui ne sont pas toutes les deux nulles en 0. Dans ce cas $(*)$ est une EDP d'inconnue $p$ qui peut se r\'e\'ecrire sous la forme $p'_x=F(x,y,p,p'_y)$ ou  $p'_y=F(x,y,p,p'_x)$ donc, dans tous les cas, le th\'eor\`eme Cauchy-Kowalevska nous dit qu'il existe une unique solution $p$ si l'on impose soit $p(0,y)$ soit $p(x,0)$. 
     Cela dit que, partant des trois fonctions d'une variable $r$, $s$ et $t$ avec $t(y)=p(0,y)$ ou $t(x)=p(x,0)$, on a une m\'ethode universelle pour construire les tissus de Nakai de rang 1. 
    Il y a cependant deux probl\`emes :

 1- il ne faut retenir parmi les $p$ trouv\'es que ceux qui ont une d\'eriv\'ee seconde par rapport \`a $x$ et $y$ non-nulle pour avoir de la courbure. 
                                                               
2- il se peut que deux donn\'ees initiales $(r,s,t)$ diff\'erentes donnent des tissus isomorphes. 
 
Nous ne savons pas comment modifier notre m\'ethode pour \'eviter ces redondances. On peut donc se poser la question : existe-t-il des tissus de Nakai de rang 1 non-isomorphes \`a l'exemple donn\'e au d\'ebut de section ? La r\'eponse est oui et on l'obtient comme suit.

On impose r constant \'egal \`a 1 et $ s(y)=(a/2)y+(1/2)y^2+(b/3)y^3.$
 Alors on a une solution $p$ de (*) avec $p(0,y)=1$. On sait calculer son 3-jet (par exemple via Maple). Elle n'est pas sous forme normale $1+xy(1+...)$, mais si on calcule le 3-jet de sa forme normale, on trouve $1+xy(1+ax+by).$
On voit ainsi que si l'on fait varier $a$ et $b$  les tissus obtenus sont deux \`a deux non-isomorphes.


\begin{thebibliography}{dango 9999}




   

 \bibitem[DL]{DL} J.P. Dufour and  D. Lehmann, {\it Rank of ordinary webs of codimension one. An effective method,} arXiv 1703.03725v1 [math.DG], 10 March 2017,  
   Pure and Applied Mathematics Quarterly, vol. XVI, nr.5, 2020.
   
   \bibitem[H]{H}   A. H\'enaut, {\it Planar web geometry through abelian relations  and connections}
   Annals of Math. 159 (2004)  425-445.
   
   \bibitem[N]{N} I. Nakai, {\it Curvature of curvilinear 4-webs and pencil of one forms}, Math. Helv. 73 (1998) , 177-205.
   
   \bibitem[Pa]{Pa}  A. Pantazi. {\it Sur la d\'etermination du rang d'un tissu plan.} C.R. Acad. Sc. Roumanie 4 (1940), 108-111.

   
   \bibitem[Pi]{Pi}   L. Pirio, {\it Equations Fonctionnelles Ab\'eliennes et G\'eom\'etrie des tissus},
  Th\`ese de doctorat de l'Universit\'e Paris VI, 2004.


\end{thebibliography}
\end{document}